# Admissible Permutations and an Algorithm of Frieze
HowardKleiman@qcc.cuny.edu

## Section 2.1 Introduction.

In this paper, we define admissible permutations and use them to reduce the running time of a 3-phase algorithm of Frieze that obtains a hamilton cycle in an extremal random directed graph from $O(n^{1.5})$ to $O(n^{\frac{4}{3}+o(1)})$ by reducing the running time of Phase 1.

## Section 2.2. Symbols, Definitions and Preliminary Theorems.

Let $\tau$ be an n-cycle in $S_n$, the symmetric group of permutations on the points, $V_n = \{1,2,...,n\}$. Then $\tau$ is a cyclic permutation. Define $T_n$ to be the set of cyclic permutations of $S_n$. Let $L_{i,j}$, $i,j \in V_n$, be independent uniform [0,1] random variables where the value of an arbitrary permutation, $\sigma$, is

$$L_\sigma = \sum_{i \in V_n} L_{i,\sigma_i}$$

Let $\sigma^* = \sigma^*(L) \in S_n$ have the property that $L_{\sigma^*} \leq L_\sigma$ for all $\sigma$ in $S_n$. Then $\sigma^*$ is an optimal solution of the assignment problem whose table of values is given by $L = [L_{i,j}]$ as $\sigma$ runs through the elements of $S_n$ and $A(L) = L_{\sigma^*}$ denotes its value. On the other hand, let $L_{\tau^*} \leq L_\tau$ for all $\tau$ in $T_n$ where $\tau$ ranges through all elements of $T_n$. Then $\tau^* = \tau^*(L)$ is an optimal solution of the random asymmetric traveling salesman problem and $T(L) = L_{\tau^*}$ denotes its value. Let $L_{i,j}$, $i,j \in V_n$ be independent uniform [0,1] random variables and consider the random variables $A_n = A(L)$ and $T_n = T(L)$. Karp [6] "patched together" the cycles of $\sigma^*$ to obtain a tour, $\tau$, of length $T_{n,1}$ such that

$$1 \leq E(A_n) \leq E(T_n) \leq E(T_{n,1}) = E(A_n) + O((\log n)^{\frac{3}{2}} n^{-0.24})$$

where $E(X)$ is the expected value of the random variable $X$.

Karp and Steele [7] simplified this algorithm by constructing a simpler $O(n^3)$ algorithm which produced a tour $\tau_2$ of length $T_{n,2}$ where

$$E(T_{n,2}) = E(A_n) + O(n^{-\frac{1}{2}}).$$

In theorem 1.1 of [1], Dyer and Frieze improved Karp and Steele's estimate to
In [4], Frieze used patching procedures from Karp [5] to devise an algorithm that obtained a hamilton circuit in a random directed graph with probability approaching 1 as $n \to \infty$. Its running time is $O(n^{1.5})$. More precisely, let $D_{n,m}$ denote the directed graph with set of vertices V = \{1,2,...,n\} and a set, E, of m arcs chosen uniformly from $K_n^D$, the complete directed graph on n vertices. Furthermore, let $m^*$ be the smallest subset of E such that the directed graph, $D_{n,m^*}$, has the property that $\delta^+(D_{n,m^*}) \geq 1$, $\delta^-(D_{n,m^*}) \geq 1$ where $\delta^+(D)$ denotes the minimum in-degree of the directed graph, D, $\delta^-(D)$ denotes the minimum out-degree of $D$. Frieze's algorithm consists of three steps. In the first one, (a), a small set of edges $E' \subseteq E_{m*}$ is constructed which almost always contains a set of about $[\log n]$ vertex disjoint cycles which cover $V$. In the second stage, (b), patching algorithms are used which allow the cycles to be "patched" into larger ones by using *2-cycle exchanges*. By the end of (b), there is almost always a cycle, $C'$, of length $n - o(n)$ plus at most $O(\log n)$ other cycles,

$$S = \{ C_i \mid i = 1, 2, ..., [c(\log n)]\}.$$

In the last step, (c), the elements of $S$ are "added" to $C'$, one by one, by a process of double rotations, to obtain a hamilton circuit. Using an algorithm of Even and Tarjan in [4], the running time of Step 1 is $O(n^{1.5})$. The running times of steps (2) and (3) are $O(n(\log n)^3)$ and $O(n^{\frac{4}{3} + o(1)})$, respectively. In this paper, the running time of Frieze's algorithm is improved by replacing his Step 1 by another procedure which obtains a set of about $\log n$ disjoint cycles covering $V$ in $O(n(\log n)^3)$ running time. Thus, the r. t. of his algorithm becomes $O(n^{\frac{4}{3} + o(1)})$ - that of Step 3 in [1]. Definitions given in the introduction apply to the remainder of the paper. From Frieze [4], we may assume in the construction of $D_{m*}$ that $m* \leq n(\log n) + kn$ where $k \to \infty$ as $n \to \infty$. More precisely, let

$$k = \begin{cases} k_1, & \text{if } k \text{ is positive,} \\ 0, & \text{if } k \text{ is zero,} \\ -k_1, & \text{if } k \text{ is negative.} \end{cases}$$

Furthermore, let

$$c = \begin{cases} e^{k_1}, & \text{if } k \text{ is positive,} \\ 1, & \text{if } k \text{ is zero,} \\ e^{-k_1}, & \text{if } k \text{ is negative.} \end{cases}$$

Then $m* \leq nl(og\, cn)$. In general, if $s \in S_n$, then $s(a)$ represents the action of the permutation $s$ on the point a. If $s \in S_n$ moves every point in $P$, then it is a *derangement*. Thus, the $n$-cycle $h$ is a derangement. An arc not lying in $D_{m*}$ is called a pseudo-arc associated with $D_{m*}$. $D_i$ $(i = 1, 2, ...)$ is a pseudo-derangement if it corresponds to a derangement, $d_i$, in $S_n$. A vertex of $D_i$ which is the initial vertex of a pseudo-arc of $D_i$ is a *pseudo-arc vertex*. All other vertices are *arc vertices*. Let $d_i$ be the derangement in $S_n$ corresponding to $D_i$. Then $s_i$ is admissible if $d_i s_i$ is also a derangement. As an example, if $d_i = (1\ 2)(3\ 4)$ and $s_i = (2\ 3)$, $d_i s_i = (1\ 2\ 3\ 4)$. Therefore, $d_i s_i = d_{i+1}$ is a derangement and thus $s_i$ is admissible. We note that in order to obtain the pseudo-derangement, $D_{i+1}$, the permutation $s_i = (2\ 3)$ is replaced by the set of arcs

$$S = \{(2,\ d_i(3)), (3,\ d_i(2))\} = \{(2,\ 4),\ (3,\ 1)\}$$

in the construction of $D_{i+1}$. In general, if $s_i$ is admissible, then $D_i s_i$ is a set of disjoint cycles covering all of the vertices in $V_C$. Since some of the arcs of $D_i$ do not necessarily lie in $D_{m*}$, $D_i s_i$ is generally a pseudo-derangement of $D_{m*}$. For simplicity, we denote running time by r. t. Given $m$ elements, using a balanced, binary search tree, we can insert or delete any element, or rebalance the tree in $O(log\ m)$ r. t.

**Lemma 2.1** $1 + \dfrac{1}{2} + \dfrac{1}{3} + ... + \dfrac{1}{n} \approx log n$

**Proof.** Using an arbitrary partition of the interval from 1 through $n$, say $(a_0, a_1, ..., a_n)$ with $1 \leq i \leq n$, the function $\dfrac{1}{i}$ is a monotonically decreasing function in each subinterval $a_{k-1} \leq i \leq a_k$. Using lower sums and upper sums, $\dfrac{1}{i}$, it follows that

$$\sum_{k=1}^{k=n} \dfrac{a_k - a_{k-1}}{a_k} \leq log(n+1) \leq \sum_{k=1}^{k=n-1} \dfrac{a_k - a_{k-1}}{a_{k-1}}$$

Therefore, if $a_i = i + 1$ for $0 \leq i \leq n$ with $a_i$ in the interval $[1, n+1]$, we obtain

$$\sum_{i=2}^{i=n+1} \dfrac{1}{i} \leq log(n+1) \leq \sum_{i=1}^{i=n} \dfrac{1}{i}$$

Since

$$log(n+1) - log(n) = log(\dfrac{n+1}{n}) = log(1 + \dfrac{1}{n}) \to log(1) = 0,$$

as $n \to \infty$, for large $n$, $\sum_{k=1}^{k=n} \dfrac{1}{k} \approx log(n)$.

**Comment.** *The best we can say about the above approximation is that it's less than 1.*

**Theorem 2.1.** (Feller [2])

*Let $s$ be a permutation randomly chosen from $S_n$. Then the number of disjoint cycles of $s$ approaches $[log\ n]$ as $n \to \infty$.*

**Proof.** Let $P = \{1, 2, ..., n, ...\}$ be the set of points of $S_n$, the symmetric group on $n$ points for large $n$. Assume that we are randomly constructing a permutation $s$ in $S_n$ such that an identity element is considered a cycle of length 1. We can describe the cycles of a given permutation on the points of $P$ by the inversions of their natural order. For instance, consider the following permutation:

$$\begin{matrix} 1 & 2 & 3 & 4 & 5 & 6 & 7 & 8 & 9 & 10 \\ 7 & 4 & 2 & 10 & 6 & 3 & 8 & 1 & 5 & 9 \end{matrix}$$

The numbers in the top row may be thought of as ordinal numbers – in this case, the natural ordering of the numbers from 1 through 10. The numbers in the bottom row are the numbers being permuted. For instance, 1 is in the eighth place, 8 is in the seventh place, 7 is in the first place. Thus, *(1 7 8)* is a disjoint cycle of the permutation. Next, 2 is in the third place, 3 is in the sixth place, 6 is in the fifth place, 5 is in the ninth place, 9 is in the tenth place, 10 is in the fourth place, 4 is in the second place. Thus, a second disjoint cycle is *(2 3 6 5 9 10 4)* is the second (and last) disjoint cycle of the permutation. Now let $s$ be a random permutation of the elements of P. At the *k-th* step of the construction of s, let $X_k = 1$ if a cycle is completed at the end of the *k*-th step; otherwise, $X_k = 0$. In general,

$$Pr(X_i = 1) = \frac{1}{n - i + 1}.$$

It follows that the number of cycles in s is

$$\Sigma_n = X_1 + X_2 + \ldots + X_n,$$

while the average number of cycles in randomly constructing a permutation is

$$m_n = 1 + \frac{1}{2} + \ldots + \frac{1}{n} \approx \log n$$

From a result of Feller in *[2]*, the number of cycles between $\log n + \alpha(\log n)^{\frac{1}{2}}$ and $\log n + \beta(\log n)^{\frac{1}{2}}$ is approximately given by $n! \{\Phi(\alpha) - \Phi(\beta)\}$. If we now assume that $\alpha = -\log(\log n)$ while $\beta = \log(\log n)$, then $\Phi(\beta) - \Phi(\alpha) \to 1$ as $n \to \infty$. Since $n^{1+\varepsilon} > n + \log n$ for any positive number $\varepsilon$ and a correspondingly large value of $n(\varepsilon) = n$, it follows that, as $n \to \infty$, the number of cycles in a randomly chosen permutation from $S_n$ approaches *[log n]*.

### Section 2.3. The Algorithm..

Before discussing the algorithm, we note that since Frieze's algorithm obtains a hamilton cycle, there can exist no cycles, $C_i = (v_1, v_2, \ldots, v_i)$ $(i = 2, 3, \ldots, m)$, $m < n$ in $D_{m*}$ such that each vertex, $v_j$, in $C_i$ has $d^-(v_i) = d^+(v_i) = 1$: *all* arcs of such cycles must lie on every hamilton cycle of length *n* which is impossible. Also, the existence of such a cycle would mean that $D_{m*}$ is not a strongly connected digraph. We use $D'_{m*}$ throughout. This allows us to include all arcs of $D_{m*}$ that must lie in *every* hamilton cycle of $D_{m*}$ in the set of cycles obtained while deleting only arcs than can lie on *no* hamilton cycle in $D_{m*}$. We construct $D_{m*}$ as a balanced, binary search tree.

(1) Before starting constructing permutations, let $d^+(x) = 1$, where *(x, y)* is the unique arc emanating from $x$. As we construct $D'_{m*}$, the contracted digraph of $D_{m*}$, *(x, y)* becomes the *2-vertex*, *xy*. If *y* also has out-degree 1, we construct a *3-vertex*, say *xyz*, etc.. Before continuing the algorithm, we delete all arcs that terminate in y. Similarly, if $d^-(v) = 1$ where the arc *(u,v)* lies in $D_{m*}$, we construct the 2-vertex *uv* and then delete all arcs that terminate in v. The set of vertices of $D'_{m*}$ is $V_C$.

(2) We next construct an *n'*-cycle, say $d_0$, in $S_n$. (We will discuss the exact construction in section 2.3.) This *n'*-cycle generally includes *2*-vertices, *3*-vertices, …, r-vertices implying that $n' < n$. If $d_0 = (a_1 a_2 \ldots a_n)$, we construct a corresponding *n'*-cycle, $D_0$, consisting of arcs or pseudo-arcs of $D'_{m*}$ that correspond to the arcs in $d_0$. i.e., if $a_1 a_2$ is an "arc" of $d_0$, then $(a_1, D_0(a_2))$ is an arc or pseudo-arc of $D'_{m*}$: If it lies in $D'_{m*}$, then it is an *arc*; otherwise, it is a *pseudo-arc*. In general, $D_0$ is a *pseudo derangement*. Some (or all) of these "arcs" may lie in $K_n^D - D_{m*}$. If $(a_1, D_0(a_2))$ is an arc of $D'_{m*}$, then $a_1$ is an *arc vertex*; otherwise, it is a *pseudo-arc vertex*. To indicate an arc, note that we use $D_0(a_2)$ rather than $d_0(a_2)$. The object of the algorithm is to transform $D_0$ into a derangement in $D'_{m*}$. i.e., a disjoint set of cycles each of whose arcs lie in $D'_{m*}$, while the sum of the sets of vertices of the cycles is $V_C$. Our way of doing this is to construct a random permutation from the arcs in $D'_{m*}$. To understand the procedure, we give a simple example. Suppose that
$$S = \{(a_1, d_0(a_i)), (a_i, d_0(a_j)), (a_j, d_0(a_1))\}$$
is a set of arcs in $D'_{m*}$ where $a_1, a_i, a_j$ are pseudo-arc vertices. Then $s = (a_1 a_i a_j)$ is a permutation in $S_n$ such that $d_1 = d_0 s$ corresponds to a new pseudo-derangement, $D_1$, of arcs where the elements of $S$ replace pseudo-arcs of $D_0$. Let $V_C = \{a_1, a_2, ..., a_{n'}\}$. Assume that $a_i < a_j$ if $i < j$. Define $D_0 = (a_1 a_2 ... a_{n'})$. For simplicity, denote $a_i$ by $i$. It follows that $D_0 = (1\ 2\ 3\ ...\ n')$. The function *ORD* represents the ordinal values of the elements of $D_0$. Thus, $ORD(1) = 1, ..., ORD(i) = i, ...$ $(1 \leq i \leq n')$. Analogously, $ORD^{-1}(a_i) = i$.

(3) We now construct the following balanced, binary search trees: *PSEUDO*, *ADD*, *DELETE*. Every pseudo-arc vertex of $D_0$ is placed in increasing order of magnitude on PSEUDO. $ORD[D_0]$ is a balanced, binary, search tree where each arc is on a branch headed by its initial vertex. $ORD^{-1}[D_0]$ is a balanced, binary search tree in which the entries are arranged in increasing order of magnitude of the domain values. Thus, we can find the $ORD^{-1}$ value of any vertex in $O(\log n)$ r. t. Let $a$ be an arbitrary element of *PSEUDO*. In the first step of the algorithm, we obtain an arc, say $[a, D_0(b)]$, lying in $D'_{m*} - D_0$. We then obtain $b$ in the following manner:

(i) $k \neq 1$.
$\quad$ If $ORD^{-1}(D_0(b)) = k$, then $ORD^{-1}(b) = k - 1$.

(ii) $k = 1$.
$\quad$ If $ORD^{-1}(D_0(b)) = 1$, then $ORD^{-1}(b) = n'$.

Finally, *ADD* is represented by two ordered sets, $ADD(i)$ and $ADD(t)$. In $ADD(i)$, the arcs of *ADD* are placed in increasing order of magnitude of their *initial* vertices; in $ADD(t)$, in increasing order of magnitude of their *terminal* vertices.

(4) We randomly choose a vertex from $PSEUDO_o$, say $a$, and randomly choose an arc from $D'_{m*} - \{D_0 \cup DELETE \cup ADD\}$, say $(a, D_0(b))$. (We define DELETE in (5).) Obtaining $ORD^{-1}(D_0(b)) = k$, if $k \neq 1$, we find $ORD^{-1}(b) = k - 1$; otherwise, $ORD^{-1}(b) = n'$. We then randomly choose arc $(b, D_0(c))$. Obtaining $ORD^{-1}(D_0(c)) = s$, let $c = ORD(s-1)$. We then obtain $(c, D_0(d))$. We continue the construction of a permutation
$$s_0 = (a b c d ... v ...)(a'b'c'd'....v'...)(a''b''...v''...)... = \prod_{j=1}^{j=m} C_j,$$

placing each corresponding arc of *D'* in both *ADD(i)* and *ADD(t)*. As we proceed, we continually check to see if the last arc chosen has an initial vertex in *ADD (i)* or a terminal vertex in *ADD(t)*. Suppose that $(u, D_0(v))$ is randomly chosen from *D'*. If an arc of form $(u, D_0(v'))$ lies in *ADD(i)*, we delete it and place in a balanced, binary search tree called *DELETE* whose elements are placed in increasing order of magnitude of the initial vertices of its elements. Similarly, suppose that $(c, D_0(v))$ lies in *ADD(t)*. Again, we delete it, and place it in *DELETE*. In the latter case, we continue the algorithm *by choosing an arc emanating from c*. We now discuss how the structure of $s_0$ evolves as the algorithm progresses. Suppose $v = d$, i.e., $(c, D_0(v)) = (c, D_0(d))$. If this occurs, we have two possibilities as shown in Figs. 2.1, 2.2. :

(a) *(v e f ... u)* is a cycle such that each of the arcs
   $(v, D_0(e))$, $(e, D_0(f))$, ... , $(u, D_0(v))$ lies in $D'_{m*}$ (Here
   $(u, D_0(v)) = (u, D_0(d))$. ) We are constructing a permutation, $s_0$, consisting of disjoint cycles. This can't be obtained if *(c, d), (u,d)* are both in it. Thus, as mentioned earlier, we delete $(c, D_0(d))$ from *ADD(i)* and *ADD(t),* add it to *DELETE*, and place $(u, D_0(v))$ in *ADD(i)* and *ADD(t)*. As we proceed in this manner, each time we place an arc in *ADD*, we check *ADD(i)* to see if the arc's initial vertex is in *PSEUDO*. If so, we delete it from *PSEUDO*. Once we have placed $(c, D_0(d))$ in *DELETE,* we continue the algorithm starting with the pseudo-arc *c*.

(b) We have already obtained at least one cycle, say $C_j$ ( $j = 1, 2, ...$).
   In this case, $v = d'$ may already belong to some cycle of $s_0$, say $C_{j'} = (d' e' f' ... r')$. Therefore, we *destroy* $C_{j'}$ by deleting the arc $(r', D_0(d'))$ from *ADD*, replacing it by $(u, D_0(v))$. Correspondingly, this yields the path *[ a,b,c,d,...,u,v,e',f',...,r']* used in the construction of $s_0$. (We note that *r* generally becomes a pseudo-arc vertex in this procedure. In that case, we would add it to *PSEUDO*.). We then randomly choose an arc of
   D' = $D'_{m*} - \{ D_0 \cup DELETE \cup ADD \}$ emanating from *c*.

   Continuing in this manner, suppose the following occurs: There exist no arcs emanationg from *c* in *D'*. Then the heading of the "*c*" branch of *D'* is changed to *cD*. We next randomly choose an arc, say $(c, D_0(z))$ from *DELETE*.

   Furthermore, if we delete an arc whose initial arc is *c,* we place it on the *cD* branch of *D'*. On the other hand, suppose that *DELETE* contains *no* arcs emanating from c. We then change the *cD* heading of *D'* to c and commence randomly choosing arcs out of D', and, as we did initially, place arcs deleted from *ADD* in *DELETE*. We continue the algorithm in this manner until PSEUDO contains no vertices, indicating that we have obtained a derangement all of whose arcs lie in $D'_{m*}$

Suppose that the number of arcs of *D'* chosen is $2(1 + \alpha) n \log(cn)$. We next use the Poisson approximation to the solution of the classical occupancy problem given in Feller, 7section *V.2, [3]*, to obtain the probability that we will run through all of the vertices in $V_C$. Let $c \to \infty$ as $n \to \infty$. If $\lambda = ne^{-\frac{(1+\alpha)n\log(cn)}{n}}$, the probability of success approaches $e^{-\lambda} = e^{-\frac{1}{c^{1+\alpha}n^\alpha}} \to e^0 = 1$, the probability of $(1+\alpha)n(\log(cn))$ random arcs passing through each vertex of $V_C$ approaches 1. Using theorem 2.1 in Feller *[3]*, it follows that the probability of obtaining a derangement, $d$, in $D'$ containing $\log n$ cycles where only one arc of $d$ is a pseudo-arc approaches 1 as $n \to \infty$. Using another $(1+\alpha)n\log(cn)$ randomly chosen arcs, we now determine the probability of obtaining a derangement containing only arcs of $D'_{m*}$. The probability that the last arc chosen in an admissible permutation, s, corresponds to an arc in $D'_{m*}$ is at least $\frac{1}{n}$. Therefore, the probability of *not* obtaining an arc is at most $1 - \frac{1}{n}$ implying that the probability of not obtaining an arc in $(1+\alpha)n\log(cn)$ is at most $(1 - \frac{1}{n})^{(1+\alpha)n\log(cn)} < e^{-(1+\alpha)\log(cn)} = \frac{1}{(cn)^{1+\alpha}}$. It follows that the probability of obtaining an arc is at least $1 - \frac{1}{(cn)^{1+\alpha}} \to 1$ as $n \to \infty$ and $c \to \infty$, concluding the proof that the probability of success of the algorithm approaches 1 as $n \to \infty$. From lemma 2.1 and theorem 2.1, as $n \to \infty$, the probability that a derangement on n' points has $[\log n']$ cycles approaches 1.

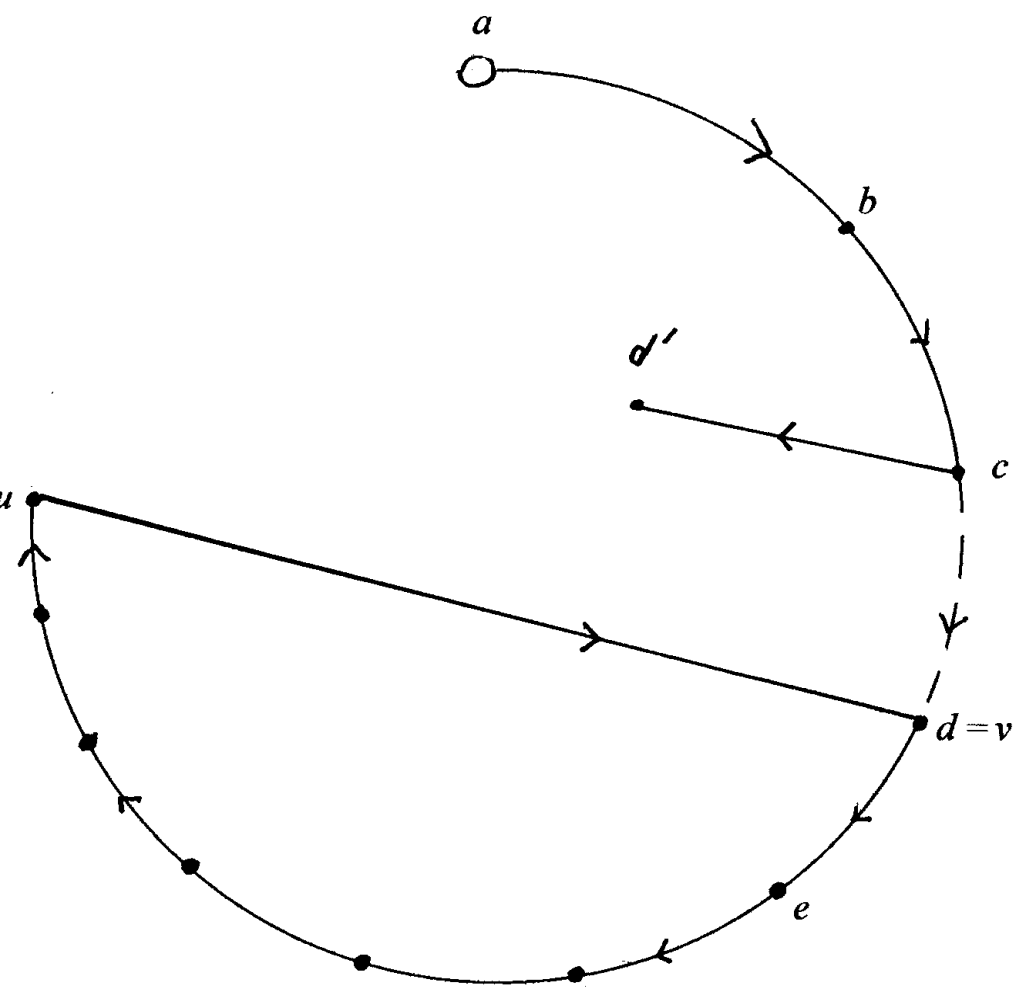

Fig. 2.1

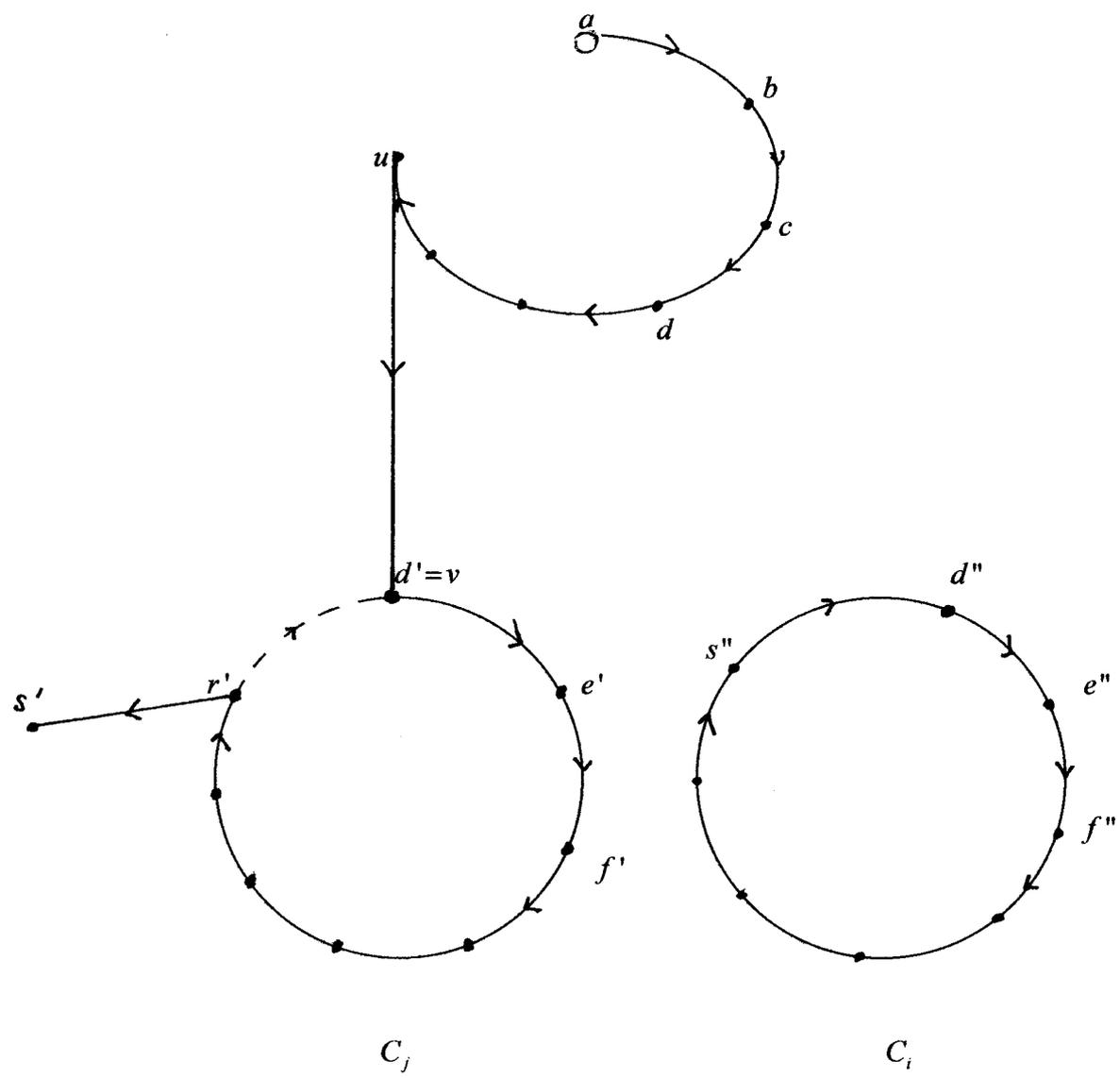

Fig. 2.2

The running time of the algorithm is $O(n(\log n)^2)$

**Example 2.1** $D'_{m*}$ is a digraph containing 27 vertices of which three are 2-vertices: 5-19, 7-25, 17-29.

$$D'_{m*}$$

1: 5-19, 13, 15

2: 9, 6

3: 2, 14, 4

4: 18, 9, 24

5-19: 23, 27, 16

6: 4, 24, 22

7-25: 8, 3, 9

8: 12, 20

9: 17-29, 26, 8

15: 6, 12, 30

16: 30, 17-29, 21

17-29: 26, 22, 18

18: 7-25, 30, 13

20: 13, 2, 27

21: 20, 18, 28

22: 3, 5-19, 2

23: 28, 1

24: 14, 7-25, 3

10: 15, 16 ,14         26: 24, 11, 17

11: 1, 10              27: 11, 1, 9

12: 10, 21, 11         28: 21, 23, 7-25

13: 27, 4, 6           30: 22, 28, 20

14: 16, 8, 10

$D_0 = $(1  2  3  4  5-19  6  7-25  8  9  10  11  12  13  14  15  16  17-29  18  20  21  22  23  24  26  27  28  30)

The following arcs of $D_0$ lie in $D'_{m^*}$: (3, 4), (7-25, 8), (16, 17-29), (17-29-18).

$ORD(1) = 1; ORD^{-1}(1) = 1$

$ORD(2) = 2; ORD^{-1}(2) = 2$

$ORD(3) = 3; ORD^{-1}(3) = 3$

$ORD(4) = 4; ORD^{-1}(4) = 4$

$ORD(5) = 5 - 19; ORD^{-1}(5 - 19) = 5$

$ORD(6) = 6; ORD^{-1}(6) = 6$

$ORD(7) = 7 - 25; ORD^{-1}(7 - 25) = 7$

$ORD(8) = 8; ORD^{-1}(8) = 8$

$ORD(9) = 9; ORD^{-1}(9) = 9$

$ORD(10) = 10; ORD^{-1}(10) = 10$

$ORD(11) = 11; ORD^{-1}(11) = 11$

$ORD(12) = 12; ORD^{-1}(12) = 12$

$ORD(13) = 13; ORD^{-1}(13) = 13$

$ORD(14) = 14; ORD^{-1}(14) = 14$

$ORD(15) = 15; ORD^{-1}(15) = 15$

$ORD(16) = 16; ORD^{-1}(16) = 16$

$ORD(17) = 17 - 29; ORD^{-1}(17 - 29) = 17$

$ORD(18) = 18; ORD^{-1}(18) = 18$

$ORD(19) = 20; ORD^{-1}(20) = 19$

$ORD(20) = 21; ORD^{-1}(21) = 20$

$ORD(21) = 22; ORD^{-1}(22) = 21$

$ORD(22) = 23; ORD^{-1}(23) = 22$

$ORD(23) = 24; ORD^{-1}(24) = 23$

$ORD(24) = 26; ORD^{-1}(26) = 24$

$ORD(25) = 27; ORD^{-1}(27) = 25$

$ORD(26) = 28; ORD^{-1}(28) = 26$

$ORD(27) = 30; ORD^{-1}(30) = 27$

$D' = D'_{m*} - \{D_0 \cup DELETE \cup ADD\}$, $DELETE = \{\}$, $ADD = \{\}$,

$PSEUDO = \{1, 2, 4, 5\text{-}19, 6, 8, 9, 10, 11, 13, 13, 14, 15, 18, 20, 21, 22, 23, 24, 26, 27, 28, 30\}$

In order to save space, we will start each random choice of an arc, (a, b), with the arc itself rather than saying " Choose (a, b)."

(1, 5-19)  1 is not in *ADD(i)*, 5-19 is not in *ADD(t)*. Place (1, 5-19) is *ADD(i)*, *ADD(t)*.
   Delete (1, 5-19) from *D'*. Delete 1 from *PSEUDO*.

$ORD^{-1}(5-19) = 5, ORD(4) = 4.$

(4, 9)   4 is not in *ADD(i)*, 9 is not in *ADD(t)*. Place (4, 9) in *ADD(i), ADD(t)*.
   Delete (4, 9) from *D'*. Delete 4 from *PSEUDO*.

$ORD^{-1}(9) = 9, ORD(8) = 8.$

(8, 12)  8 is not in *ADD(i)*, 12 is not in *ADD(t)*. Place (8, 12) in *ADD(i), ADD(t)*. ).
   Delete (8, 12) from *D'*. Delete 8 from PSEUDO.

$ORD^{-1}(12) = 12, ORD(11) = 11.$

(11, 10)  11 is not in *ADD(i)*, 10 is not in *ADD(t)*. Place (11,10) in *ADD(i), ADD(t)*.
   Delete (11, 10) from *D'*.

$ORD^{-1}(10) = 10, ORD(9) = 9.$

(9, 17-29)  9 is not in *ADD(i)*, 17-29 is not in *ADD(t)*. Place (9, 17-29) in *ADD(i)*, *ADD(t)*. Delete (9, 17-29) from *D'*. Delete 9 from *PSEUDO*.

PSEUDO = {2, 5-19, 6, 10, 12, 13, 14, 15, 18, 20, 21, 22, 23, 24, 26, 27, 28, 30}

*ADD(i)* = {(1, 5-19), (4, 9), (8, 12), (9, 17,29), (11, 10)}

*ADD(t)* = {(1, 5-19), (4, 9), (11,10), (8, 12), (9, 17-29)}

ARCS IN $D_0$ = {(3, 4), (7-25, 8), (16, 17-29), (17-29, 18)}

DELETE = { }

$ORD^{-1}(17-29) = 17, ORD(16) = 16.$

(16, 21)  16 is not in *ADD(i)*, 21 is not in *ADD(t)*. Place (16, 21) in
   *ADD(i), ADD(t)*.
   Delete (16, 21) from D'.

$ORD^{-1}(21) = 20, ORD(19) = 20.$

(20, 13)  20 is not in *ADD(i)*, 13 is not in *ADD(t)*. Place (20, 13) in
   *ADD(i), ADD(t)*.
   Delete (20, 13) from *D'*. Delete 20 from *PSEUDO*.

$ORD^{-1}(13) = 13, ORD(12) = 12.$

(12, 21)  12 is not in *ADD(i)*, 21 is in *ADD(t)*: (16, 21). Delete (16, 21) from
   *ADD(i), ADD(t)*. Place (12, 21) in *ADD(i), ADD(t)*. Delete 12 from
   *PSEUDO*.

(16, 30)  Place (16, 21) in *DELETE*. 16 is not in *ADD(i)*, 30 is not in *ADD(t)*. Place
   (16, 30) in *ADD(i), ADD(t)*. Delete (16, 30) from *D'*.

$ORD^{-1}(30) = 27, ORD(26) = 28.$

(28, 7-25)  28 is not in *ADD(i)*, 7-25 is not in *ADD(t)*. Place (28, 7-25) in *ADD(i)*,
        *ADD(t)*. Delete (28, 7-25) from *D'*. Delete 28 from *PSEUDO*.

$ORD^{-1}(7-25)=7, ORD(6)=6.$

(6, 4)      6 is not in *ADD(i)*, 4 is not in *ADD(t)*. Place (6, 4) in *ADD(i), ADD(t)*.
        Delete it from *D'*. Delete 6 from *PSEUDO*.

$$PSEUDO = \{2, 5\text{-}19, 10, 13, 14, 15, 18, 21, 22, 23, 24, 26, 27, 30\}$$

$$ADD(i) = \{(1, 5\text{-}19), (4, 9), (6, 4), (8, 12), (9, 17\text{-}29), (11, 10), (12, 21), (16, 30),\\ (20, 13), (28, 7\text{-}25)\}$$

$$ADD(t) = \{(6, 4), (1, 5\text{-}19), (28, 7\text{-}25), (4, 9), (11, 10), (8, 12), (20, 13), (9, 17\text{-}29),\\ (12, 21), (16, 30)\}$$

$$DELETE = \{(6, 4), (4, 9), (16, 21)\}$$

$$ARCS \text{ in } D_0 = \{(3, 4), (7\text{-}25, 8), (16, 17\text{-}29), (17\text{-}29, 18)\}$$

$ORD^{-1}(4)=4, ORD(3)=3.$

(3, 14)    3 doesn't lie in *ADD(i)*, 14 doesn't lie in *ADD(t)*. Place (3, 14) in
        *ADD(i), ADD(t)*. Delete (3, 14) from *D'*.

$ORD^{-1}(14)=14, ORD(13)=13.$

(13, 4)    13 is not in *ADD(i)*, 4 is in *ADD(t)*: (6, 4). Delete (6, 4) from
        *ADD(i), ADD(t)*.
        Place (13, 4) in *ADD(i), ADD(t)*. Delete (13, 4) from *D'*.

(6, 24)    Place (6, 4) in *DELETE*. 6 is not in *ADD(i)*, 24 is not in *ADD(t)*. Place
        (6, 24) in *ADD(i), ADD(t)*. Delete it from *D'*.

$ORD^{-1}(24)=23, ORD(22)=23.$

(23, 28)  23 is not in *ADD(i)*, 28 is not in *ADD(t)*. Place (23, 28) in *ADD(i), ADD(t)*.
        Delete it from *D'*. Delete 23 from *PSEUDO*.

$ORD^{-1}(28)=26, ORD(25)=27.$

(27, 29)  27 is not in *ADD(i)*, 9 is in *ADD(t)*: (4, 9). Delete (4, 9) from
        *ADD(i), ADD(t)*.
        Place (27, 9) in *ADD(i), ADD(t)*. Delete it from *D'*. Delete 27 from *PSEUDO*.

$$PSEUDO = \{2, 4, 5\text{-}19, 10, 14, 15, 18, 21, 22, 24, 26, 30\}$$

$$ADD(i) = \{(1, 5\text{-}19), (3, 14), (6, 24), (8, 12), (9, 17\text{-}29), (11. 10), (12, 21), (13, 4),\\ (16, 30), (20, 13), (23, 28), (27, 9), (28, 7\text{-}25)\}$$

$$ADD(t) = \{(13, 4), (1, 5\text{-}19), (28, 7\text{-}25), (27, 9), (11, 10), (8, 12), (20, 13), (3, 14),\\ (9, 17\text{-}29), (12, 21), (23, 28), (28, 7\text{-}25), (16, 30)\}$$

$$DELETE = \{(4, 9), (6, 4), (16, 21)\}$$

$$ARCS \ in \ D_0 \ = \ \{(3, 4), (7\text{-}25, 8), (16, 17\text{-}29), (17\text{-}29, 18)\}$$

(4, 24)   Place (4, 9) in *DELETE*. 4 is not in *ADD(i)*, 24 is in *ADD(t)*: (6, 24). Delete (6, 24) from *ADD(i), ADD(t)*. Place (4, 24) in *ADD(i), ADD(t)*. Delete (4, 24) from *D'*. Delete 4 from *PSEUDO*.

(6, 22)   Place (6, 24) in *DELETE*. 6 is not in *ADD(i)*, 22 is not in *ADD(t)*. Place (6, 22) in *ADD(i), ADD(t)*. Delete (6, 22) from *D'*.

$ORD^{-1}(22) = 21, ORD(20) = 21.$

(21, 20) 21 is not in *ADD(i)*, 20 is not in *ADD(t)*. Place (21, 20) in *ADD(i), ADD(t)*. Delete (21, 20) from *D'*. Delete 21 from *PSEUDO*.

$ORD^{-1}(20) = 19, ORD(18) = 18.$

(18, 30) 18 is not in *ADD(i)*, 30 is in *ADD(t)*: (16, 30). Delete (16, 30) from *ADD(i), ADD(t)*. Place (18, 30) in ADD(i), ADD(t). Delete (18, 30) from *D'*. Delete 18 from *PSEUDO*.

No arcs emanate from 16 in *D'*. Rename the heading of the branch 16 in *D'*, 16D. Choose an arc emanating from 16 in *DELETE*.

(16, 21)   16 is not in *ADD(i)*, 21 is in *ADD(t)*: (12, 21). Delete (12, 21) from *ADD(i), ADD(t)*. Place (16, 21) in *ADD(i), ADD(t)*.

(12, 11)   Place (12, 21) in *DELETE*. 12 is not in *ADD(i)*, 11 is not in *ADD(t)*. Place (12, 11) in *ADD(i), ADD(t)*. Delete (12, 11) from *D'*

$ORD^{-1}(11) = 11, ORD(10) = 10..$

$$PSEUDO = \{2, 5\text{-}19, 10, 14, 15, 22, 24, 26, 30\}$$

$$ADD(i) = \{(1, 5\text{-}19), (3, 14), (4, 24), (6, 22), (8, 12), (9, 17\text{-}29), (11, 10), (12, 11),$$
$$(13, 4), (16, 21), (18, 30), (20, 13), (21, 20), (23, 28), (27, 9), (28, 7\text{-}25)\}$$

$$ADD(t) = \{(13, 4), (1, 5\text{-}19), (28, 7\text{-}25), (27, 9), (11, 10), (12, 11), (8, 12), (20, 13),$$
$$(3, 14), (9, 17\text{-}29), (21, 20), (16, 21), (6, 22), (4, 24), (23, 28), (18, 30)\}$$

$$DELETE = \{(3, 4), (4, 9), (6, 4), (6, 24), (12, 21), (16, 30)\}$$

$$ARCS \ in \ D_0 \ = \ \{(3, 4), 7\text{-}25, 8), (16, 17\text{-}29), (17\text{-}29, 18)\}$$

(10, 16)   10 is not in *ADD(i)*, 16 is not in *ADD(t)*. Place (10, 16) in *ADD(i), ADD(t)*. Delete (10, 16) from *D'*. Delete 10 from *PSEUDO*.

$ORD^{-1}(16) = 16, ORD(15) = 15.$

(15, 12)   15 is not in *ADD(i)*, 12 is in ADD(t): (8, 12). Delete (8, 12) from *ADD(i), ADD(t)*. Place (15, 12) in *ADD(i), ADD(t)*. Delete (15, 12) from

    *D'.* Delete 15 from *PSEUDO.*
(8, 20)  8 is not in *ADD(i)*, 20 is in *ADD(t)*: (21, 20). Delete (21, 20) from
    *ADD(i), ADD(t).* Place (8, 20) in *ADD(i), ADD(t).* Delete (8, 20) from *D'.*
(21, 18) Place (21, 20) in *DELETE.* 21 is not in *ADD(i),* 18 is not in *ADD(t).* Place
    (21, 18) in *ADD(i), ADD(t).* Delete (21, 18) from *D'.*
$ORD^{-1}(18)=18, ORD(17)=17-29.$
(17-29, 26) 17-29 is not in *ADD(i),* 26 is not in *ADD(t).* Place (17-29, 26) in
    *ADD(i), ADD(t).* Delete (17-29, 26) from *D'.*
$ORD^{-1}(26)=24, ORD(23)=24.$

$$PSEUDO = \{2, 5\text{-}19, 14, 22, 24, 26, 30\}$$

$$\begin{aligned}ADD(i) = \{&(1, 5\text{-}19), (3, 14), (4, 24), (6, 22), (8, 20), (9, 17\text{-}29), (10, 16), (11, 10),\\&(12, 11), (15, 12), (16, 21), (17\text{-}29, 26), (18, 30), (20, 13), (21, 18),\\&(23, 28), (27, 9), (28, 7\text{-}25)\}\end{aligned}$$
$$\begin{aligned}ADD(t) = \{&(13, 4), (1, 5\text{-}19), (28, 7\text{-}25), (27, 9), (11, 10), (12, 11), (15, 12), (20, 13),\\&(3, 14), (10, 16), (9, 17\text{-}29), (21, 18), (8, 20), (16, 21),\\&(6, 22), (4, 24), (17\text{-}29, 26), (23, 28), (18, 30)\}\end{aligned}$$

$$DELETE = \{(4, 9), (6, 4), (6, 24), (8, 12), (16, 30), (21, 20)\}$$

$$ARCS\ in\ D_0 = \{(3, 4), (7\text{-}25, 8), (16, 17\text{-}29), (17\text{-}29, 18)\}$$

(24, 3) 24 is not in *ADD(i),* 3 is not in *ADD(t).* Place (24, 3) in *ADD(i), ADD(t).*
   Delete (24, 3) from *D'* Delete 24 from *PSEUDO..*
$ORD^{-1}(3)=3, ORD(2)=2.$
(2, 6)  2 is not in *ADD(i),* 6 is not in *ADD(t).* Place (2, 6) in *ADD(i), ADD(t).*
   Delete (2, 6) from *D'.* Delete 2 from *PSEUDO.*
$ORD^{-1}(6)=6, ORD(5)=5-19.$
(5-19, 23) 5-19 is not in *ADD(i),* 23 is not in *ADD(t).* Place (5-19, 23) in
   *ADD(i), ADD(t).* Delete (5-19, 23) from D'. Delete 5-19 from *PSEUDO.*
$ORD^{-1}(23)=22, ORD(21)=22.$
(22, 2) 22 is not in *ADD(i),* 2 is not in *ADD(t).* Place (22, 2) in *ADD(i), ADD(t).*
   Delete (22, 2) from *D'.* Delete 22 from *PSEUDO.*
$ORD^{-1}(2)=2, ORD(1)=1.$
(1, 13) 1 is in *ADD(i)*: (1, 5-19), 13 is in *ADD(t)*: (20, 13). Delete (1, 5-19), (20, 13)
   from *ADD(i), ADD(t).* Place (1, 13) in *ADD(i), ADD(t).* Delete (1, 13) from
   *D'.*
(20, 27) Place (1, 5-19), (20, 13) in *DELETE.* 20 is not in *ADD(i),* 27 is not in
   *ADD(t).* Place (20, 27) in *ADD(i), ADD(t).* Delete (20, 27) from *D'..*
$ORD^{-1}(27)=25, ORD(24)=26.$

$$PSEUDO = \{14, 26, 30\}$$

$ADD(i) = \{(1, 13), (2, 6), (3, 14), (4, 24), (5\text{-}19, 23), (6, 22), (8, 20), (9, 17\text{-}29),$
$(10, 16), (11, 10), (12, 11), (13, 4), (15, 12), (16, 21), (17\text{-}29, 26),$
$(18, 30), (20, 27), (21, 18), (22, 2), (23, 28), (24, 3), (27, 9), (28, 7\text{-}25)\}$

$ADD(t) = \{(22, 2), (24, 3), (13, 4), (2, 6), (28, 7\text{-}25), (27, 9), (11, 10), (12, 11),$
$(15, 12), (1, 13), (3, 14), (10, 16), (9, 17\text{-}29), (21, 18), (8, 20), (16, 21),$
$(6, 22), (5\text{-}19, 23), (4, 24), (20, 27), (23, 28), 18, 30)\}$

$DELETE = \{(1, 5\text{-}19), (4, 9), (6, 4), (8, 12), (16, 30), (20, 13), (21, 20)\}$

$ARCS\ in\ D_0 = \{(3, 4), (7\text{-}25, 8), 16, 17\text{-}29), (17\text{-}29, 18)\}$

(26, 24)  26 is not in *ADD(i)*, 24 is in *ADD(t)*: (4, 24). Delete (4, 24) from *ADD(i), ADD(t)*. Place (26, 24) in ADD(i), ADD(t). Delete (26, 24) from *D'*. Delete 26 from *PSEUDO*.

(4, 18)  Place (4, 24) in *DELETE*. 4 is not in *ADD(i)*, 18 is in *ADD(t)*: (21, 18). Delete (21, 18) from *ADD(i), ADD(t)*. Place (4, 18) in *ADD(i), ADD(t)*. Delete (4, 18) from *D'*.

(21, 28)  Place (21, 18) in *DELETE*. 21 is not in *ADD(i)*, 28 is in *ADD(t)*: (23, 28). Delete (23, 28) from *ADD(i), ADD(t)*. Place (21, 28) in *ADD(i), ADD(t)*. Delete (21, 28) from *D'*.

(23, 1)  Place (23, 28) in *DELETE*. 23 is not in *ADD(i)*, 1 is not in *ADD(t)*. Place (23, 1) in *ADD(i), ADD(t)*. Delete (23, 1) from *D'*.

$ORD^{-1}(1)=1, ORD(27)=30.$

(30, 20)  30 is not in *ADD(i)*, 20 is in *ADD(t)*: (8, 20). Delete (8, 20) from *ADD(i), ADD(t)*. Place (30, 20) in *ADD(i), ADD(t)*. Delete (30, 20) from *D'*. Delete 30 from *PSEUDO*.

No arcs emanate out of 8 in *D'*. Rename the heading 8 of *D'*, 8D.
Choose (8, 12) from *DELETE*. Delete (8, 12) from *DELETE*. Place (8, 20) in *D'*.
  8 is not in *ADD(i)*, 12 is in *ADD(t)*: (15, 12). Delete (15, 12) from *ADD(i), ADD(t)*. Place (8, 12) in *ADD(i), ADD(t)*. Place 15 in *PSEUDO*.

$$PSEUDO = \{14, 15\}$$

$ADD(i) = \{(1, 13), (2, 6), (3, 14), (4, 18), (5\text{-}19, 23), (6, 22), (8, 12), (9, 17\text{-}29),$
$(10, 16), (11, 10), (12, 11), (13, 4), (16, 21), (17\text{-}29, 26), (18, 30), (20, 27),$
$(21, 28), (22, 2), (23, 1), (24, 3), (26, 24), (27, 9), (28, 7\text{-}25), (30, 20)\}$

$ADD(t) = \{((23, 1), (22, 2), (24, 3), (13, 4), (2, 6), (28, 7\text{-}25), (27, 9), (11, 10), (12, 11),$
$(8, 12), (1, 13), (3, 14), (10, 16), (9, 17\text{-}29), (4, 18), (30, 20), (16, 21),$

(6, 22), (5-19, 23), (26, 24), (17-29, 26), (20, 27), (21, 28), (18, 30)}

*DELETE* = {(1, 5-19), (4, 9), (4, 24), (6, 4), (6, 24), (8, 12), (16, 30), (20, 13), (21, 18)}

*ARCS in* $D_0$ = {(3, 4), (7-25, 8), (16, 17-29), (17-29, 18)}

(15, 6)   Place (15, 12) in *DELETE*. 15 is not in *ADD(i)*, 6 is in *ADD(t)*: (2, 6). Delete
(2, 6) from *ADD(i), ADD(t)*. Place (15, 6) in *ADD(i), ADD(t)*. Delete
(15, 6) from *D'*. Delete 15 from *PSEUDO*.

(2, 9)    Place (2, 6) in *DELETE*. 2 is not *ADD(i)*, 9 is in *ADD(t)*: (27, 9). Delete
(27, 9) from *ADD(i), ADD(t)*. Place (2, 9) in *ADD(i), ADD(t)*. Delete (2, 9)
from *D'*.

(27, 1)   Place (27, 9) in *DELETE*. 27 is not in *ADD(i)*, 1 is in *ADD(t)*: (23, 1). Delete
(23, 1) from *ADD(i), ADD(t)*. Place (27, 1) in *ADD(i), ADD(t)*. Delete
(27, 1) from *D'*.

No arcs emanate from 23 in *D'*. Rename the heading 23 in D', 23D.
Choose (23, 28) from *DELETE*. Place (23, 1) in *D'*.
      23 is not in *ADD(i)*, 28 is in *ADD(t)*: (21, 28). Delete (21, 28) from
      *ADD(i), ADD(t)*. Place (23, 28) in *ADD(i), ADD(t)*.

(21, 20)  Place (21, 28) in *DELETE*. 21 is not in *ADD(i)*, 20 is in *ADD(t)*: (30, 20).
Delete (30, 20) from *ADD(i), ADD(t)*. Place (21, 20) in *ADD(i), ADD(t)*..
Delete (21, 20) from *D'*. Place 30 in *PSEUDO*.

*PSEUDO* = {14, *30*}

*ADD(i)* = {(1, 13), (2, 9), (3, 14), (4, 18), (5-19, 23), (6, 220, (8, 12), (9, 17-29),
    (10, 16), (11, 10), (12, 11), (13, 4), (15, 6), (16, 21), (17-29, 26), (18, 30),
    (20, 27), (21, 20), (22, 2), (23, 28), (24, 3), (26, 24), (27, 1), (28, 7-25)}

*ADD(t)* = {(27, 1), (22, 2), (24, 3), (13, 4), (15, 6), (28, 7-25), (2, 9), (11, 10),
    (12, 11), (8, 12), (1, 13), (3, 14), (10, 16), (9, 17-29), (4, 18), (21, 20),
    (16, 21), (6, 22), (5-19, 23), (26, 24), (17-29, 26), (20, 27), (23, 28), (18, 30)}

*DELETE* = {(1, 5-19), (2, 6), (4, 9), 4, 24), (6, 4), (6, 24), (16, 30), (20, 13), (21, 18),
                                                    (21, 28), (23, 28)}

*ARCS in* $D_0$ = {(3, 4), (7-25, 8), (16, 17-29), (17-29, 18)}

(30, 22)   Place (30, 20) in *DELETE*. 30 is not in *ADD(i)*, 22 is in *ADD(t)*: (6, 22).
       Delete (6, 22) from *ADD(i), ADD(t)*. Place (30, 22) in *ADD(i), ADD(t)*.
       Delete 30 from *PSEUDO*.

No arc emanates from 6 in *D'*. Rename 6, 6D in *D'*.
Choose (6, 4) from *DELETE*. Place (6, 22) in *D'*.

|  | 6 is not in *ADD(i)*, 4 is in *ADD(t)*: (13, 4). Delete (13, 4) from *ADD(i), ADD(t)*. Place (6, 4) in *ADD(i), ADD(t)*. Delete (6, 4) from *DELETE*. |
|---|---|
| (13, 6) | Place (13, 4) in *DELETE*. 13 is not in *ADD(i)*, 6 is in *ADD(t)*: (15, 6). Delete (15, 6) from *ADD(i), ADD(t)*. Place (13, 6) in *ADD(i), ADD(t)*.. |
| (15, 30) | Place (15, 6) in *DELETE*. 15 is not in *ADD(i)*, 30 is in *ADD(t)*: (18, 30). Delete (18, 30) from *ADD(i), ADD(t)*. Place (15, 30) in *ADD(i), ADD(t)*.. |
| (18, 13) | Place (18, 30) in *DELETE*. 18 is not in *ADD(i)*, 13 is in *ADD(t)*: (1, 13). Delete (1, 13) from *ADD(i), ADD(t)*. Place (18, 13) in *ADD(i), ADD(t)*.. |
| (1, 15) | Place (1, 13) in *DELETE*. 1 is not in *ADD(i)*, 15 is not in *ADD(t)*. Place (1, 15) in *ADD(i), ADD(t)*. |

$ORD^{-1}(15) = 15$, $ORD(14) = 14$.

| (14, 8) | 14 is not in *ADD(i)*, 8 is not in *ADD(t)*. Place (14, 8) in *ADD(i), ADD(t)*. Delete 14 from *PSEUDO*. |
|---|---|

$ORD^{-1}(8) = 8$, $ORD(7) = 7 - 25$.

| (7-25, 3) | 7-25 is not in *ADD(i)*, 3 is in *ADD(t)*: (24, 3). Delete (24, 3) from *ADD(i), ADD(t)*. Place (7-25, 3) in *ADD(i), ADD(t)*. Place 24 in *PSEUDO*. |
|---|---|

$$PSEUDO = \{24\}$$

$ADD(i) = \{(1, 15), (2, 9), (3, 14), (4, 18), (5\text{-}19, 23), (6, 4), (7\text{-}25, 3), (8, 12),$
$(9, 17\text{-}29), (10, 16), (11, 10), (12, 11), (13, 6), (14, 8), (15, 30),$
$(16, 21), (17\text{-}29, 26), (18, 13), (20, 27), (21, 20), (22, 2),$
$(23, 28), (26, 24), (27, 1), (28, 7\text{-}25), (30, 22)\}$

$ADD(t) = \{(27, 1), (22, 2), (7\text{-}25, 3), (6, 4), (13, 6), (28, 7\text{-}25), (14, 8), (2, 9),$
$(11, 10), (12, 11), (8, 12), (18, 13), (3, 14), (1, 15), (10, 16),$
$(9, 17\text{-}29), (4, 18), (21, 20), (16, 21), (30, 22), (5\text{-}19, 23),$ .
$(26, 24), (27, 1), (28, 7\text{-}25), (30, 22)\}$

$DELETE = \{(1, 5\text{-}19), (1, 13), (2, 6), (4, 9), (4, 24), (6, 4), (6, 24), (8, 12),$
$(13, 4), (15, 6), (16, 30), (18, 30), (20, 13), (21, 18), (21, 28),$
$(23, 28), (24, 3), 30, 20)\}$

$$ARCS \text{ in } D_0 = \{(3, 4), (7\text{-}25, 8), (16, 17\text{-}29), (17\text{-}29, 18)\}$$

| (24, 14) | Place (24, 3) in *DELETE*. 24 is not in *ADD(i)*, 14 is in *ADD(t)*: (3, 14). Delete (3, 14) from *ADD(i), ADD(t)*. Place (24, 14) in *ADD(i), ADD(t)*. Delete 24 from *PSEUDO*. |
|---|---|
| (3, 2) | Place (3, 14) in *DELETE*. 3 is not in *ADD(i)*, 2 is in *ADD(t)*: (22, 2). |

              Delete (22, 2) from *ADD(i), ADD(t)*. Place (3, 2) in *ADD(i), ADD(t)*.
(22, 5-19)  Place (22, 2) in *DELETE*. 22 is not in *ADD(i)*, 5-19 is not in *ADD(t)*.

We thus have no more vertices in *PSEUDO*.

$$PSEUDO = \{\}$$

*ADD(i)* = {(1, 15), (2, 9), (3, 2), (4, 18), (5-19, 23), (6, 4), (7-25, 3), (8, 12),
        (9, 17-29), (10, 16), (11, 10), (12, 11), (13, 6), (14, 8), (15, 30),
        (16, 21), (17-29, 26), (18, 13), (20, 27), (21, 20), (22, 5-19),
        (23, 28), (24, 14), (26, 24), (27, 1), (28, 7-25), (30, 22)}

D = (1  15  30  22  5-19  23  28  7-25  3  2  9  17-29  26  24  14  8  12  11  10  16  21  20  27)
            (4  18  13  6)

.